\documentclass[12pt]{article}

\usepackage{amssymb,amsmath,amsthm,graphicx}

\newcommand{\alf}{\alpha}

\newcommand{\del}{\delta}
\newcommand{\eps}{\varepsilon}

\newcommand{\tet}{\theta}

\newcommand{\lam}{\lambda}

\newcommand{\abs}[1]{\left|#1\right|}

\newcommand{\lang}{\big\langle}
\newcommand{\rang}{\big\rangle}

\newcommand{\bN}{\mathbb{N}}
\newcommand{\bR}{\mathbb{R}}

\newcommand{\bT}{\mathbb{T}}
\newcommand{\bZ}{\mathbb{Z}}

\begin{document}

\begin{center}
{\large Finite searches, Chowla's cosine problem, \\
and large Newman polynomials}
\end{center}

\begin{flushright}
Idris Mercer\\
Florida International University \\
\verb+imercer@fiu.edu+
\end{flushright}

\begin{abstract}
A length $n$ cosine sum is an expression of the form
$\cos a_1\tet + \cdots + \cos a_n\tet$ where $a_1 < \cdots < a_n$
are positive integers,
and a length $n$ Newman polynomial is an expression of the form
$z^{a_1} + \cdots + z^{a_n}$ where $a_1 < \cdots < a_n$
are nonnegative integers.
We define
$-\lam(n)$ to be the largest minimum of a length $n$ cosine sum
as $\{a_1,\ldots,a_n\}$ ranges over all sets of $n$ positive integers,
and we define
$\mu(n)$ to be the largest minimum modulus on the unit circle
of a length $n$ Newman polynomial as $\{a_1,\ldots,a_n\}$
ranges over all sets of $n$ nonnegative integers.
Since there are infinitely many possibilities for the $a_j$, it is
not obvious how to compute $\lam(n)$ or $\mu(n)$
for a given $n$ in finitely many steps.
Campbell et al.\ found the value of $\mu(3)$ in 1983,
and Goddard found the value of $\mu(4)$ in 1992.
In this paper, we find the values of $\lam(2)$ and $\lam(3)$
and nontrivial bounds on $\mu(5)$.
We also include further remarks on the seemingly difficult
general task of reducing the computation of $\lam(n)$ or $\mu(n)$
to a finite problem.
\end{abstract}

\section{Introduction}

We define a {\bf length $n$ cosine sum} to be any expression of the form
$$
\cos a_1\tet + \cdots + \cos a_n\tet
$$
where $a_1 < \cdots < a_n$ are positive integers.
If $f(\tet)$ is a length $n$ cosine sum, then $f(0)=n$ and $f(\tet)\ge-n$ for all $\tet$.
For $f(\tet)=\cos a_1\tet + \cdots + \cos a_n\tet$, we define
$$
L(a_1,\ldots,a_n) = \min_\tet f(\tet) = \min_{0\le\tet\le\pi} f(\tet)
$$
which is negative, since $f(\tet)$ has average value 0 on $[0,\pi]$.
More specifically, we have
$$
-n \le L(a_1,\ldots,a_n) < 0.
$$
We then define
$$
-\lam(n) = \sup L(a_1,\ldots,a_n)
$$
where $\{a_1,\ldots,a_n\}$ ranges over all of the infinitely many sets of
$n$ positive integers. Then $\lam(n)$ is a well-defined function of $n$
because we are taking the supremum of a bounded set. Specifically, we have
$0 \le \lam(n) \le n$.
There is also a construction showing $\lam(n)=O(\sqrt{n})$
(see Section 2).
However, since there are infinitely many possibilities
for the $a_j$, it is not obvious how to compute $\lam(n)$ for a given $n$
in finitely many steps.

We define a {\bf length $n$ Newman polynomial} to be any expression of the form
$$
z^{a_1} + \cdots + z^{a_n}
$$
where $a_1 < \cdots < a_n$ are nonnegative integers.
If $f(z)$ is a length $n$ Newman polynomial and $\abs{z}=1$,
then $0 \le \abs{f(z)} \le n$. We define
$$
M(a_1,\ldots,a_n) = \min_{\abs{z}=1} \abs{z^{a_1}+\cdots+z^{a_n}}
$$
so we have
$$
0 \le M(a_1,\ldots,a_n) \le n
$$
and in fact one can show $M(a_1,\ldots,a_n)\le\sqrt{n}$ by considering the $L^2$ norm of~$f$.
(We note that the case $n=2$ is uninteresting because
$z^{a_1}+z^{a_2}=z^{a_1}(1+z^{a_2-a_1})$ always has minimum modulus 0 at the
$(a_2-a_1)$th roots of $-1$.)
We then define
$$
\mu(n) = \sup M(a_1,\ldots,a_n)
$$
where $\{a_1,\ldots,a_n\}$ ranges over all of the infinitely many sets of
$n$ nonnegative integers. Then $\mu(n)$ is a well-defined function of $n$
because we are taking the supremum of a bounded set.
However, since there are infinitely many possibilities
for the $a_j$, it is not obvious how to compute $\mu(n)$ for a given $n$
in finitely many steps.

One can inquire about the growth rates of the functions $\lam(n)$ and $\mu(n)$,
or particular values of $\lam(n)$ and $\mu(n)$.

Around the late 1940s, Ankeny and Chowla conjectured \cite{Cho52,Cho65}
that $\lam(n)$ approaches infinity with $n$.
This was first proved by Uchiyama and Uchiyama \cite{UU}
using results of Cohen \cite{Coh};
their lower bound for $\lam(n)$ was sublogarithmic.
Over the years, better lower bounds for $\lam(n)$ have been found.
The best lower bound currently known is due to Ruzsa \cite{Ruz};
it is superlogarithmic but grows more slowly than any power of~$n$.
The best known upper bound for $\lam(n)$ appears to be $O(\sqrt{n})$.
Chowla conjectured \cite{Cho65} that this is the true rate of growth.

It seems that the growth of $\mu(n)$ is less studied than the growth of $\lam(n)$.
There is a construction showing $\mu(n)$ exceeds a power of $n$
for infinitely many $n$; specifically, we have $\mu(n) \ge n^{0.14}$
when $n$ is a power of 9. See Section 2.
Boyd \cite{Boy} conjectured that $\mu(n)>1$ for all $n\ge6$
and that $\log\mu(n)/\log n$ approaches a positive constant
as $n$ approaches infinity.
The current author \cite{Mer} proved that $\mu(n)>0$ for all $n>2$.
It appears that nobody has proved that $\mu(n)$ approaches infinity
with $n$.

Although there are infinitely many cosine sums and Newman polynomials
of a given length, some brute-force exploration of specific examples
can lead to conjectures about particular values of $\lam(n)$
and $\mu(n)$.

Specifically, some experimentation leads to the conjectures
\begin{align*}
\lam(2) & = -L(1,2) = 9/8 = 1.125, \\
\lam(3) & = -L(1,2,3) \approx 1.315565, \\
\lam(4) & = -L(1,2,3,4) \approx 1.519558, \\
\lam(5) & = -L(1,2,4,5,6) \approx 1.627461, \\
\lam(6) & = -L(1,2,4,6,7,8) \approx 1.591832,
\end{align*}
and
\begin{align*}
\mu(3) & = M(0,1,3) \approx 0.607346, \\
\mu(4) & = M(0,1,2,4) \approx 0.752394, \\
\mu(5) & = M(0,1,2,6,9) = 1.
\end{align*}
Campbell et al.\ \cite{CFF} proved that
$\mu(3)=M(0,1,3)$, and Goddard \cite{God} proved that
$\mu(4)=M(0,1,2,4)$.
In this article, we prove that $\lam(2)=-L(1,2)$
and $\lam(3)=-L(1,2,3)$, and that
$1 \le \mu(5) \le 1 + \pi/5$.
To prove results of this type, one must reduce a
potentially infinite search to a finite search.

It would be interesting to show that for any $n$,
the value of $\lam(n)$ or $\mu(n)$ can be computed in a
finite number of steps (even a ridiculously large
finite number). It appears to be difficult to prove this.

\section{Constructions bounding $\lam(n)$ and $\mu(n)$}

There are some straightforward constructions that lead to
bounds on $\lam(n)$ and $\mu(n)$. These probably count as
mathematical folklore. We include them here for completeness.

{\bf Fact 2.1.} For each positive integer $n$, there exists a
length $n$ cosine sum $f(\tet)$ that satisfies
$f(\tet) \ge -\sqrt{2n}-\frac12$ for all $\tet$.
It follows that $\lam(n) \le \sqrt{2n}+\frac12$.

{\it Proof.}
Let $k=\lceil\sqrt{2n}\rceil$, so $\sqrt{2n}\le k<\sqrt{2n}+1$.
Let $\{b_1<\ldots<b_k\}$ be a set of $k$ nonnegative integers
such that the $k \choose 2$ positive differences $b_j-b_i$
are all distinct; for instance, we can take $b_j = 2^{j-1}$.
(A set with this property is sometimes called a `Sidon set';
in this proof, we do not need our Sidon set to be optimal in
any sense.) Then if $z=e^{i\tet}$, we have
\begin{align*}
0 \le \abs{z^{b_1}+\cdots+z^{b_k}}^2 & =
(z^{b_1}+\cdots+z^{b_k})(z^{-b_1}+\cdots+z^{-b_k}) \\
& = k + \sum_{i<j} 2\mathrm{Re}(z^{b_j-b_i})
= k + 2\sum_{i<j}\cos(b_j-b_i)\tet.
\end{align*}
Thus if we define $g(\tet)=\sum_{i<j}\cos(b_j-b_i)\tet$,
then $g(\tet)$ is a length $k \choose 2$ cosine sum
satisfying $g(\tet)\ge-k/2$ for all $\tet$.
Now observe
$$
n - \sqrt{\frac{n}{2}} =
\frac{\sqrt{2n}(\sqrt{2n}-1)}{2} \le
\frac{k(k-1)}{2} < \frac{(\sqrt{2n}+1)\sqrt{2n}}{2}
= n + \sqrt{\frac{n}{2}}.
$$
If ${k \choose 2} > n$, choose $h(\tet)$ to be a sum
consisting of ${k \choose 2} - n$ of the cosines
in the sum $g(\tet)$, and choose $f(\tet)=g(\tet)-h(\tet)$.
If ${k \choose 2} < n$, choose $h(\tet)$ to be a sum
of $n - {k \choose 2}$ cosines not appearing in the sum $g(\tet)$,
and choose $f(\tet) = g(\tet) + h(\tet)$.
If ${k \choose 2} = n$, define $h(\tet)=0$ and choose $f(\tet)=g(\tet)$.

Then $h(\tet)$ is a sum of at most $\sqrt{n/2}$ cosines,
so we have $\pm h(\tet) \ge -\sqrt{n/2}$. We then have
$$
f(\tet) = g(\tet)\pm h(\tet) \ge -\frac{k}{2}-\sqrt{\frac{n}{2}}
> -\frac{\sqrt{2n}+1}{2}-\sqrt{\frac{n}{2}}
= -\sqrt{2n}-\frac12
$$
which completes the proof.

{\bf Fact 2.2.} Let $f$ be a length $n$ Newman polynomial satisfying $\abs{f(z)}\ge K_1$
for all $\abs{z}=1$, and let $g$ be a length $m$ Newman polynomial satisfying $\abs{g(z)}\ge K_2$
for all $\abs{z}=1$. Then there exists a length $nm$ Newman polynomial $h$ satisfying
$\abs{h(z)}\ge K_1K_2$ for all $\abs{z}=1$.

{\it Proof.} Suppose
\begin{align*}
f(z) & = z^{a_1} + \cdots + z^{a_n}, \\
g(z) & = z^{b_1} + \cdots + z^{b_m}.
\end{align*}
If $k$ is a sufficiently large positive integer, then the product
$$
h(z) = f(z^k)g(z) = (z^{ka_1}+\cdots+z^{ka_n})(z^{b_1}+\cdots+z^{b_m})
$$
has the property that the $nm$ exponents $ka_i+b_j$ are all distinct,
and is hence a Newman polynomial of length $nm$.
If $\abs{z}=1$, then $\abs{h(z)}=\abs{f(z^k)}\abs{g(x)}\ge K_1K_2$.
This completes the proof.

As mentioned in \cite{Boy}, the length 9 Newman polynomial
$$
f(z) = 1 + z + z^2 + z^3 + z^4 + z^7 + z^8 + z^{10} + z^{12}
$$
has unusually high minimum modulus on the unit circle.
Specifically, $\abs{f(z)} \ge 1.362$ for all $z$ on the unit circle.
Then, by repeated application of Fact 2.2, we can construct for each $k$
a Newman polynomial $f(z)$ of length $9^k$ with the property that
$\abs{f(z)}\ge1.362^k$ for all $z$ on the unit circle.
Since $1.362>9^{0.14}$, this means that for infinitely many
values of $n$, we have a length $n$ Newman polynomial $f(z)$
that satisfies $\abs{f(z)}\ge n^{0.14}$ on the unit circle.

For example, there is a Newman polynomial of length $9^3=729$
that satisfies $\abs{f(z)}\ge1.362^3\approx2.53$ on the unit circle,
so $\mu(729)\ge2.53$. It would be interesting to know
(for example) the least $n$ such that $\mu(n)\ge2$.

\section{Some notation and lemmas}

In this section, we establish some notation and some useful lemmas.

Let $f(\tet)=\cos a_1\tet + \cdots + \cos a_n\tet$,
let $g=\gcd(a_1,\ldots,a_n)$, and let $a'_j = a_j/g$.
Then $\cos a'_1\tet + \cdots + \cos a'_n\tet$ is a
length $n$ cosine sum taking on the same values
(and hence having the same minimum) as $f(\tet)$.
Therefore, in the definition of $\lam(n)$, it suffices to consider
only those $\{a_1,\ldots,a_n\}$ for which $\gcd(a_1,\ldots,a_n)=1$.

Let $f(z) = z^{a_1} + \cdots + z^{a_n}$, where $a_1<\cdots<a_n$
are nonnegative integers. Define $b_j=a_j-a_1$. If $\abs{z}=1$,
then
$$
f_1(z) = f(z)/z^{a_1} = 1 + z^{b_2} + \cdots + z^{b_n}
$$
has the same modulus as $f(z)$. Next, define
$g=\gcd(b_2,\ldots,b_n)$ and $b'_j = b_j/g$.
If
$$
f_2(z) = 1 + z^{b'_2} + \cdots + z^{b'_n}
$$
then $f_2(z^g)=f_1(z)$, so $f_1(z)$ and $f_2(z)$ have the same
set of outputs as $z$ ranges over the unit circle.
Therefore, in the definition of $\mu(n)$, it suffices to consider
only those $\{a_1,\ldots,a_n\}$ of the form $\{0,a_2,\ldots,a_n\}$
where $\gcd(a_2,\ldots,a_n)=1$.
Furthermore, there is one more symmetry we exploit.
If
$$
f(z) = 1 + z^{a_2} + \cdots + z^{a_{n-1}} + z^{a_n}
$$
then we define
$$
h(z) = z^{a_n}f(z^{-1})
= z^{a_n} + z^{a_n-a_2} + \cdots + z^{a_n-a_{n-1}} + 1
$$
and note that $\abs{f(z)}$ and $\abs{h(z)}$ have the same set of
outputs as $z$ ranges over the unit circle.
So we are free to choose between $f(z)$ and $h(z)$ and can hence
assume $a_{n-1}\ge a_n-a_2$.

For the above reasons, we make the following definitions.
Define
\begin{align*}
\bN'_n & = \{(a_1,\ldots,a_n) \mid 0<a_1<\cdots<a_n, \;
\gcd(a_1,\ldots,a_n)=1\}, \\
\bN''_n & = \{(0,a_2,\ldots,a_n) \mid 0<a_2<\cdots<a_n, \;
\gcd(a_2,\ldots,a_n)=1, \; a_{n-1}\ge a_n-a_2\}
\end{align*}
and note that if $L(a_1,\ldots,a_n)$ and $M(a_1,\ldots,a_n)$
are as defined in Section 1,
then we have
$$
-\lam(n) = \sup L(a_1,\ldots,a_n)
$$
where the supremum is taken over all $(a_1,\ldots,a_n)\in\bN'_n$,
and we have
$$
\mu(n) = \sup M(0,a_2,\ldots,a_n)
$$
where the supremum is taken over all $(0,a_2,\ldots,a_n)\in\bN''_n$.

We also define $\bT$ to be $\bR$ mod $2\pi$.
Following are some definitions and lemmas regarding subsets of $\bT$.

{\bf Definition.}
An {\bf equispaced subset} of $\bT$ of {\bf order} $m$ is
any subset of $\bT$ of the form
$$
\Big\{ \tet_0 + \frac{2k\pi}{m} \;\Big|\; k\in\bZ \Big\}.
$$
Note that if we fix $\xi\in\bT$, then the set $\{\tet\in\bT\mid m\tet=\xi\}$
is an equispaced set of order $m$.

{\bf Lemma 3.1.}
Let $S_1$ and $S_2$ be equispaced subsets of $\bT$ of order $m_1$ and $m_2$
respectively, and let $g=\gcd(m_1,m_2)$. Then there exists $\tet_1\in S_1$
and $\tet_2\in S_2$ such that $\abs{\tet_1-\tet_2}\le\frac{\pi g}{m_1m_2}$.

{\it Proof.}
Suppose
$$
S_1 = \Big\{ \xi_1 + \frac{2k\pi}{m_1} \;\Big|\; k\in\bZ \Big\} \qquad \mbox{and}
\qquad S_2 = \Big\{ \xi_2 + \frac{2k\pi}{m_2} \;\Big|\; k\in\bZ \Big\}.
$$
The real number $\frac{m_1m_2}{2\pi}(\xi_1-\xi_2)$ must be within $\frac{g}{2}$ of an
integer multiple of $g$, call it $ag$.
Also, $ag$ can be written in the form $km_1-\ell m_2$ for some integers $k$ and $\ell$.
Thus we have
$$
\abs{ \frac{m_1m_2}{2\pi}(\xi_1-\xi_2) - (km_1-\ell m_2) } \le \frac{g}{2}
$$
which, multiplying by $\frac{2\pi}{m_1m_2}$ and rearranging, gives us
$$
\abs{ \Big(\xi_1+\frac{2\pi\ell}{m_1}\Big)
- \Big(\xi_2+\frac{2\pi k}{m_2}\Big) } \le \frac{\pi g}{m_1m_2},
$$
completing the proof of the lemma.

The following lemma is straightforward.

{\bf Lemma 3.2.}
If $\tet$ satisfies $\abs{\tet-\pi}=\eps$, then $\cos\tet \le -1+\frac{1}{2}\eps^2$.

We will also need other bounds on the cosine function,
which need not be the best bounds possible.
One can show the following.

{\bf Lemma 3.3.}
If $\tet$ satisfies $\abs{\tet-\frac{2\pi}{3}}=\eps\le\frac{\pi}{6}$,
then
$$
\cos\tet \le -\frac12 +\frac{3}{\pi}\eps
$$
and if $\tet$ satisfies $\abs{\tet-\frac{4\pi}{3}}=\eps\le\frac{\pi}{6}$,
then
$$
\cos\tet \le -\frac12 +\frac{3}{\pi}\eps.
$$

{\bf Definition.}
If $S\subset\bT$ is an equispaced set of order $m$, and $f(\tet)$ and $g(\tet)$ are real-valued
functions on $\bT$, then we define
$$
\lang f(\tet),g(\tet)\rang_S = \frac1m\sum_{\tet\in S}f(\tet)g(\tet),
$$
the average value of $f(\tet)g(\tet)$ over $S$, which we can think of as a kind of dot product
of $f(\tet)$ and $g(\tet)$.

One can verify that this dot product has the following properties:
\begin{itemize}
\item $\lang f_1(\tet)+f_2(\tet),g(\tet)\rang_S
= \lang f_1(\tet),g(\tet)\rang_S + \lang f_2(\tet),g(\tet)\rang_S$
\item $\lang f(\tet),g_1(\tet)+g_2(\tet)\rang_S
= \lang f(\tet),g_1(\tet)\rang_S + \lang f(\tet),g_2(\tet)\rang_S$
\item $\lang1,1\rang_S = 1$
\item $\lang1,\cos k\tet\rang_S = \lang\cos k\tet,1\rang_S = 0$ if $k$ is {\bf not} a multiple of $m$
\item $\lang\cos k\tet,\cos k\tet\rang_S = \frac12$ if $2k$ is {\bf not} a multiple of $m$
\item $\lang\cos k\tet,\cos\ell\tet\rang_S = 0$ if $k+\ell$ and $k-\ell$ are {\bf not} multiples of $m$
\end{itemize}

A function that is nonnegative on $\bT$ can be used as a `weight function'.
Some examples of nonnegative weight functions are:
\begin{gather*}
1 - \cos k\tet \\
2(1 - \cos k\tet)^2
= 3 - 4\cos k\tet + \cos2k\tet
\end{gather*}
as well as any sum of such functions.

{\bf Fact 3.4.}
Let $w(\tet)$ be a nonnegative weight function on $\bT$,
let $g(\tet)$ be any real-valued function on $\bT$,
and let $S\subset\bT$ be an equispaced set.
If we have $\lang w(\tet),g(\tet)\rang_S\le0$, then $g(\tet)\le0$
for some $\tet\in\bT$.

\section{The values of $\lam(2)$ and $\lam(3)$}

To prove $\lam(2) = -L(1,2) = 9/8$, let $(a,b)\in\bN'_2$, so $\gcd(a,b)=1$.
We must show that $\cos a\tet + \cos b\tet \le -9/8$ for some $\tet$.

If $b\le2$, then $\cos a\tet + \cos b\tet = \cos\tet + \cos2\tet$,
which one can verify has minimum value $-9/8$. So suppose $b\ge3$.
Define
$$
S_1 = \{ \tet\in\bT \mid a\tet=\pi \} \qquad \mbox{and} \qquad S_2 = \{ \tet\in\bT \mid b\tet=\pi \}
$$
which are equispaced sets of order $a$ and $b$ respectively.
By Lemma 3.1, there exist $\tet_1\in S_1$ and $\tet_2\in S_2$
such that $\abs{\tet_1-\tet_2}\le\frac{\pi}{ab}$.
Let $\tet=\tet_2$. Then $b\tet=\pi$, and
$\abs{a\tet-\pi}=\abs{a\tet_2-a\tet_1}\le\frac{\pi}{b}$,
so by Lemma 3.2, we have
$$
\cos a\tet \le -1 + \frac{\pi^2}{2b^2} \le -1 + \frac{\pi^2}{18} \approx -0.45, 
$$
implying $\cos a\tet+\cos b\tet \le -2 + \frac{\pi^2}{18} \approx -1.45 < -9/8$.
This completes the evaluation of $\lam(2)$.

Next, we will show that $\lam(3) = -L(1,2,3)$.
First, observe that trigonometric identities allow us to write
\begin{align*}
\cos\tet + \cos2\tet + \cos3\tet & = \cos\tet + (2\cos^2\tet-1) + (4\cos^3\tet-3\cos\tet) \\
& = 4\cos^3\tet + 2\cos^2\tet - 2\cos\tet - 1
\end{align*}
which is a polynomial in $\cos\tet$ of degree 3. One can verify that
$$
\min_{-1\le c\le1} 4c^3 + 2c^2 - 2c - 1 = -\frac{17+7\sqrt{7}}{27} \approx -1.315565.
$$
That is, $-L(1,2,3)=\frac{17+7\sqrt{7}}{27}\approx1.315565.$
For brevity, let $K=-L(1,2,3)$.

To prove that $\lam(3)=K$, we will partition $\bN'_3$ into four subsets:
\begin{align*}
M_1 & = \{ (a,b,c)\in\bN'_3 \mid c=2a \}, \\
M_2 & = \{ (a,b,c)\in\bN'_3 \mid c=2b \}, \\
M_3 & = \{ (a,b,c)\in\bN'_3 \mid c=a+b \}, \\
M_0 & = \big\{ (a,b,c)\in\bN'_3 \mid c \notin \{2a,2b,a+b\} \big\}.
\end{align*}
Speaking very informally, we can think of $M_1,M_2,M_3$ as subsets of $\bN'_3$
that have only two `degrees of freedom'.
For each $j\in\{0,1,2,3\}$, we will show that with finitely many exceptions,
if $(a,b,c)\in M_j$ then $\cos a\tet+\cos b\tet+\cos c\tet\le-K$ for some $\tet$.

Suppose $(a,b,c)\in M_0$. We then have $0<a<b<c$, $\gcd(a,b,c)=1$, and $c\notin\{2a,2b,a+b\}$.
Let $S=\{\tet\in\bT\mid c\tet=\pi\}$, which is an equispaced set of order $c$.
Note that $2-\cos a\tet-\cos b\tet$ is a nonnegative weight function, and consider
$$
\del = \Big\langle 2-\cos a\tet-\cos b\tet, \; \frac12 + \cos a\tet + \cos b\tet \Big\rangle_S.
$$
Now observe the following:
\begin{itemize}
\item Since $0<a<c$, we conclude $a$ is not a multiple of $c$
\item Since $0<b<c$, we conclude $b$ is not a multiple of $c$
\item Since $0<b-a<c$, we conclude $b-a$ is not a multiple of $c$
\item Since $0<2a<2c$, if $2a$ is a multiple of $c$ then $2a=c$
\item Since $0<2b<2c$, if $2b$ is a multiple of $c$ then $2b=c$
\item Since $0<a+b<2c$, if $a+b$ is a multiple of $c$ then $a+b=c$
\end{itemize}
It then follows from properties of the dot product that
$$
\del = \lang2,\frac12\rang_S - \lang\cos a\tet,\cos a\tet\rang_S - \lang\cos b\tet,\cos b\tet\rang_S
= 1 - \frac12 - \frac12 = 0
$$
and then by Fact 3.4, we conclude that $\frac12+\cos a\tet+\cos b\tet \le 0$
for some $\tet\in S$.
But that $\tet$ satisfies $\cos c\tet=-1$ and therefore $\cos a\tet+\cos b\tet+\cos c\tet \le -3/2 < -K$.

Next, suppose $(a,b,c)\in M_1$. If $a\le2$ then $(a,b,c)$ must be $(2,3,4)$,
and we observe that $\cos2\tet+\cos3\tet+\cos4\tet=-2<-K$ if $\tet=\pi/3$.
We therefore assume $a\ge3$.
Now define
$$
S_1 = \Big\{\tet\in\bT \mid a\tet=\frac{2\pi}{3}\Big\}, \qquad S_2 = \{\tet\in\bT \mid b\tet=\pi\},
$$
which are equispaced sets of order $a$ and $b$ respectively.
Note that $\gcd(a,b)=1$ since any common divisor would divide $2a=c$.
By Lemma 3.1, there exist $\tet_1\in S_1$ and $\tet_2\in S_2$ such that
$\abs{\tet_1-\tet_2}\le\frac{\pi}{ab}$.
Let $\tet=\tet_1$. Then $\cos a\tet=\cos2\pi/3=-1/2$ and $\cos c\tet=\cos2a\tet=\cos4\pi/3=-1/2$.
Also, we have $\abs{b\tet-\pi}=\abs{b\tet_1-b\tet_2}\le\frac{\pi}{a}$,
so by Lemma 3.2, we have
$$
\cos b\tet \le -1+\frac{\pi^2}{2a^2} \le -1+\frac{\pi^2}{18} \approx -0.45,
$$
implying $\cos a\tet + \cos b\tet + \cos c\tet \le -2+\frac{\pi^2}{18} \approx -1.45 < -K$.

Next, suppose $(a,b,c)\in M_2$. This case is very similar to $(a,b,c)\in M_1$.
If $b\le2$ then $(a,b,c)$ must be $(1,2,4)$, and we observe that
$\cos\tet+\cos2\tet+\cos4\tet=-3/2<-K$ if $\tet=2\pi/3$. We therefore assume $b\ge3$.
We define
$$
S_1 = \{\tet\in\bT \mid a\tet=\pi\}, \qquad S_2 = \Big\{\tet\in\bT \mid b\tet=\frac{2\pi}{3}\Big\},
$$
which are equispaced sets of order $a$ and $b$, and we note that $\gcd(a,b)=1$,
so there exist $\tet_1\in S_1$ and $\tet_2\in S_2$ such that $\abs{\tet_1-\tet_2}\le\frac{\pi}{ab}$.
We choose $\tet=\tet_2$. Then $\cos b\tet=\cos2\pi/3=-1/2$ and $\cos c\tet=\cos 2b\tet=\cos4\pi/3=-1/2$.
Also, we have $\abs{a\tet-\pi}=\abs{a\tet_2-a\tet_1}\le\frac{\pi}{b}$, so Lemma 3.2 gives us
$$
\cos a\tet \le -1 + \frac{\pi^2}{2b^2} \le -1 + \frac{\pi^2}{18} \approx -0.45,
$$
implying $\cos a\tet+\cos b\tet+\cos c\tet\le-2+\frac{\pi^2}{18}\approx-1.45<-K$.

Finally, suppose $(a,b,c)\in M_3$. We will show the existence of a finite subset $M' \subset M_3$
such that if $(a,b,c)\in M_3\setminus M'$, then $\cos a\tet+\cos b\tet+\cos c\tet\le-K$ for some $\tet$.
That will reduce this case to checking the finitely many elements of $M'$.
It will suffice to choose
$$
M' = \{(a,b,c)\in M_3 \mid b\le32\} = \{ (a,b,a+b)\in M_3 \mid 1\le a<b\le32 \}
$$
which has at most ${32 \choose 2}=496$ elements. At the end of this section, we will give an alternative argument
that avoids such a large finite set.

So suppose $(a,b,c)\in M_3\setminus M'$, so $c=a+b$ and $b\ge33$. We choose
$$
S_1 = \Big\{\tet\in\bT \mid a\tet=\frac{2\pi}{3}\Big\}, \qquad
S_2 = \Big\{\tet\in\bT \mid b\tet=\frac{2\pi}{3}\Big\},
$$
which are equispaced sets of order $a$ and $b$ respectively.
Note that $\gcd(a,b)=1$ since any common divisor would divide $a+b=c$.
By Lemma 3.1, there exist $\tet_1\in S_1$ and $\tet_2\in S_2$
such that $\abs{\tet_1-\tet_2}\le\frac{\pi}{ab}$.
Let $\tet=\tet_2$. Then $\cos b\tet=\cos b\tet_2=\cos2\pi/3=-1/2$.
Also, we have $\abs{a\tet-2\pi/3}=\abs{a\tet_2-a\tet_1}\le\frac{\pi}{b}$
and
$$
\abs{c\tet-\frac{4\pi}{3}}=\abs{\Big(a\tet-\frac{2\pi}{3}\Big)+\Big(b\tet-\frac{2\pi}{3}\Big)}
= \abs{\Big(a\tet-\frac{2\pi}{3}\Big)+0} \le \frac{\pi}{b}.
$$
By Lemma 3.3, we then have
$$
\cos a\tet \le -\frac12 + \frac3\pi\frac\pi{b} \qquad \mbox{and} \qquad
\cos c\tet \le -\frac12 + \frac3\pi\frac\pi{b}
$$
which implies
$$
\cos a\tet+\cos b\tet+\cos c\tet \le -\frac32 + \frac6b \le -\frac32 + \frac6{33} \approx -1.318 < -K.
$$

We remark that the case $(a,b,c)\in M_3$ can be dealt with more easily
if we take as known the fact that
$$
\mu(3) = M(0,1,3) = \sqrt{\frac{47-14\sqrt{7}}{27}} \approx 0.607346,
$$
which was shown in Section 3 of \cite{CFF}.
Note that this result says that given any degree 3 Newman polynomial
$f(z)=1+z^{k}+z^{\ell}$, we have $\abs{f(z)}\le M(0,1,3)$ for some $\abs{z}=1$.
We apply this to the case $f(z)=1+z^a+z^{a+b}$.
The statement
$$
\abs{1+z^a+z^{a+b}} \le \sqrt{\frac{47-14\sqrt{7}}{27}} \qquad \mbox{for some $\abs{z}=1$}
$$
is equivalent to
$$
\abs{1+z^a+z^{a+b}}^2 \le \frac{47-14\sqrt{7}}{27} \qquad \mbox{for some $\abs{z}=1$}
$$
but we also have, if $z=e^{i\tet}$,
\begin{align*}
\abs{1+z^a+z^{a+b}}^2 & = (1+z^a+z^{a+b})(1+z^{-a}+z^{-a-b}) \\
& = 3 + (z^a+z^{-a}) + (z^b+z^{-b}) + (z^{a+b}+z^{-a-b}) \\
& = 3 + 2\big(\cos a\tet + \cos b\tet + \cos(a+b)\tet \big)
\end{align*}
and the condition
$$
3 + 2\big(\cos a\tet + \cos b\tet + \cos(a+b)\tet \big) \le \frac{47-14\sqrt{7}}{27}
$$
is equivalent to
$$
\cos a\tet + \cos b\tet + \cos(a+b)\tet \le \frac12\Big(\frac{47-14\sqrt{7}}{27}-3\Big)
= -\frac{17+7\sqrt{7}}{27} = L(1,2,3).
$$
In \cite{CFF}, the fact that $\mu(3)=M(0,1,3)$ is proved by reducing to a finite number
of cases, but their number of cases is much less than 496.
Since we used the value of $\mu(3)$ just in the case where $(a,b,c)\in M_3$,
it could be said that our evaluation of $\lam(3)$ is a generalization of
the evaluation of $\mu(3)$.

\section{A possible strategy for $\lam(4)$}

The current author is unaware of how to evaluate $\lam(4)$, but includes in this section
a possible outline of a strategy where we reduce that problem to a finite set of problems
that, speaking informally, have fewer `degrees of freedom'.

We conjecture that $\lam(4) = -L(1,2,3,4) \approx 1.519558$.
To prove this, we must show that if $(a,b,c,d)\in\bN'_4$,
then $\cos a\tet+\cos b\tet+\cos c\tet+\cos d\tet \le L(1,2,3,4)$ for some $\tet$.

Let $(a,b,c,d)\in\bN'_4$, and define $S=\{\tet\in\bT \mid d\tet=\pi\}$,
which is an equispaced set of order $d$.
Define the nonnegative weight function
\begin{align*}
w(\tet) & = (1-\cos a\tet) + (1-\cos b\tet) + 2(1-\cos c\tet)^2 \\
& = 1-\cos a\tet + 1-\cos b\tet + 3 - 4\cos c\tet + \cos2c\tet \\
& = 5 - \cos a\tet - \cos b\tet - 4\cos c\tet + \cos2c\tet.
\end{align*}
If we can show that
$$
\Big\langle w(\tet), \; \frac35 + \cos a\tet + \cos b\tet + \cos c\tet \Big\rangle_S \le 0
$$
then it will follow that $\cos a\tet+\cos b\tet+\cos c\tet \le -3/5$
for some $\tet$ satisfying $\cos d\tet=-1$, so then $\cos a\tet+\cos b\tet+\cos c\tet+\cos d\tet\le-8/5=-1.6<L(1,2,3,4)$.

Using previously stated properties of the dot product, we note that evaluating
$$
\Big\langle 5 - \cos a\tet - \cos b\tet - 4\cos c\tet + \cos2c\tet, \; \frac35 + \cos a\tet + \cos b\tet + \cos c\tet \Big\rangle_S
$$
depends on properties of the set $\{0,a,b,c,2c\}\pm\{0,a,b,c\}$.
We can illustrate using the following table.
\begin{center}
\begin{tabular}{c|ccccc}
 & 0 & $a$ & $b$ & $c$ & $2c$ \\ \hline
0 & 0 & $a$ & $b$ & $c$ & $2c$ \\
$a$ & $a$ & $2a$ & $b\pm a$ & $c\pm a$ & $2c\pm a$ \\
$b$ & $b$ & $b\pm a$ & $2b$ & $c\pm b$ & $2c\pm b$ \\
$c$ & $c$ & $c\pm a$ & $c\pm b$ & $2c$ & $2c\pm c$ \\
\end{tabular}
\end{center}
If all the nonzero numbers in the body of the table are nonmultiples of $d$,
then properties of the dot product give us
\begin{align*}
& \Big\langle 5 - \cos a\tet - \cos b\tet - 4\cos c\tet + \cos2c\tet, \; \frac35 + \cos a\tet + \cos b\tet + \cos c\tet \Big\rangle_S \\
= \; & 5\cdot\frac35 - \frac12 - \frac12 - 4\cdot\frac12 = 0
\end{align*}
and so in that case, the desired conclusion follows. We now observe:
\begin{itemize}
\item The positive numbers $a,\;b,\;c,\;b-a,\;c-a,\;c-b$ are all less than $d$, and are hence nonmultiples of $d$.
\item The positive numbers $2a,\;2b,\;2c,\;b+a,\;c+a,\;c+b,\;2c-a,\;2c-b$ are all less than $2d$, and hence if any
of them are multiples of $d$, they are equal to $d$.
\item The positive numbers $2c+a,\;2c+b,\;2c+c$ are all less than $3d$, and hence if any of them are multiples of $d$,
they are equal to $d$ or $2d$.
\end{itemize}
We thus have a finite list of 14 conditions
\begin{align*}
d & = 2a \\
d & = 2b \\
d & = 2c \\
d & = b+a \\
d & = c+a \\
d & = c+b \\
d & = 2c-a \\
d & = 2c-b \\
d & = 2c+a \\
2d & = 2c+a \\
d & = 2c+b \\
2d & = 2c+b \\
d & = 3c \\
2d & = 3c
\end{align*}
such that if $(a,b,c,d)\in\bN'_4$ satisfies none of these 14 conditions,
then $\cos a\tet+\cos b\tet+\cos c\tet+\cos d\tet\le-8/5<L(1,2,3,4)$ for some $\tet$.
It then remains to deal with those $(a,b,c,d)$ that do satisfy some of these
14 conditions. If we define
\begin{align*}
M_1 & = \{(a,b,c,d)\in\bN'_4 \mid d = 2a\} \\
M_2 & = \{(a,b,c,d)\in\bN'_4 \mid d = 2b\} \\
M_3 & = \{(a,b,c,d)\in\bN'_4 \mid d = 2c\} \\
M_4 & = \{(a,b,c,d)\in\bN'_4 \mid d = b+a\} \\
M_5 & = \{(a,b,c,d)\in\bN'_4 \mid d = c+a\} \\
M_6 & = \{(a,b,c,d)\in\bN'_4 \mid d = c+b\} \\
M_7 & = \{(a,b,c,d)\in\bN'_4 \mid d = 2c-a\} \\
M_8 & = \{(a,b,c,d)\in\bN'_4 \mid d = 2c-b\} \\
M_9 & = \{(a,b,c,d)\in\bN'_4 \mid d = 2c+a\} \\
M_{10} & = \{(a,b,c,d)\in\bN'_4 \mid d = c+\frac12a\} \\
M_{11} & = \{(a,b,c,d)\in\bN'_4 \mid d = 2c+b\} \\
M_{12} & = \{(a,b,c,d)\in\bN'_4 \mid d = c+\frac12b\} \\
M_{13} & = \{(a,b,c,d)\in\bN'_4 \mid d = 3c\} \\
M_{14} & = \{(a,b,c,d)\in\bN'_4 \mid d = \frac32c\}
\end{align*}
then each $M_j$, speaking very informally, is a subset of $\bN'_4$
having only 3 `degrees of freedom'.
If we can prove for each $j$ that if $(a,b,c,d)\in M_j$,
then $\cos a\tet+\cos b\tet+\cos c\tet+\cos d\tet\le L(1,2,3,4)$ for some $\tet$,
then that would complete the evaluation of $\lam(4)$.
Perhaps it is possible to reduce each $M_j$ to a finite collection
of subsets of $\bN'_4$ that have 2 degrees of freedom.

\section{Bounds on $\mu(5)$}

It was observed in \cite{God} that
$$
M(0,1,2,6,9) = \min_{\abs{z}=1}\abs{1+z+z^2+z^6+z^9} = 1
$$
and it is suspected that $\mu(5)=1$, although this appears to be
difficult to prove.
In this section, we give a short elementary argument
showing that $\mu(5)\le\sqrt{3}\approx1.732$,
and a longer case analysis showing that
$\mu(5)\le1+\pi/5\approx1.628$.
More generally, we will show that for any positive integer $m$,
the problem of showing $\mu(5)\le1+\pi/m$ can be
reduced to checking a finite number of cases.

To show $\mu(5)\le K$, we must show that for all of the
infinitely many $(0,a,b,c,d)$ in $\bN''_5$, we have
$\abs{1+z^a+z^b+z^c+z^d}\le K$ for some $\abs{z}=1$.

Let $(0,a,b,c,d)\in\bN''_5$.
To show $\abs{1+z^a+z^b+z^c+z^d}\le\sqrt{3}$ for some $\abs{z}=1$,
define $S=\{\tet\in\bT \mid d\tet=\pi\}$,
and let $z=e^{i\tet}$.
Observe that we have
\begin{align*}
\abs{z^a+z^b+z^c}^2 & = (z^a+z^b+z^c)(z^{-a}+z^{-b}+z^{-c}) \\
& = 3 + (z^{b-a}+z^{a-b}) + (z^{c-a}+z^{a-c}) + (z^{c-b}+z^{b-c}) \\
& = 3 + \cos(b-a)\tet + \cos(c-a)\tet + \cos(c-b)\tet.
\end{align*}
Now note that $b-a,\;c-a,\;c-b$ are positive integers less than $d$,
so they are not multiples of $d$.
It follows that the average of $\abs{z^a+z^b+z^c}^2$ over $S$
is $3$.
Therefore at least one $\tet\in S$ satisfies
$\abs{z^a+z^b+z^c}^2\le3$ and hence
$\abs{1+z^a+z^b+z^c+z^d}=\abs{z^a+z^b+z^c}\le\sqrt{3}$.

To get better bounds on $\mu(5)$, we will use some
straightforward lemmas similar to those in Section 3.

{\bf Lemma 6.1.}
If $\tet$ satisfies $\abs{\tet-\pi}\le\del$,
then $\abs{1+e^{i\tet}}\le\del$.

{\bf Lemma 6.2.}
Let $k,\ell,m$ be distinct positive integers with $\ell<m$,
and let $g=\gcd(k,m-\ell)$. We then have
$$
\abs{1+z^k+z^\ell+z^m} \le \frac{\pi g}{k}
$$
for some $z$ on the unit circle.

{\it Proof.}
Define the sets
\begin{align*}
S_1 & = \{\tet\in\bT \mid k\tet=\pi\}, \\
S_2 & = \{\tet\in\bT \mid (m-\ell)\tet=\pi\},
\end{align*}
which are equispaced sets of order $k$ and $m-\ell$
respectively.
By Lemma 3.1, there exist $\tet_1\in S_1$ and $\tet_2\in S_2$
such that $\abs{\tet_1-\tet_2}\le\frac{\pi g}{k(m-\ell)}$, and hence
$$
\abs{(m-\ell)\tet_1-\pi} = \abs{(m-\ell)\tet_1-(m-\ell)\tet_2}
\le\frac{\pi g}{k}.
$$
Then, if $z=e^{i\tet_1}$, we have
\begin{align*}
\abs{1+z^k+z^\ell+z^m} & = \abs{1+z^k+z^\ell(1+z^{m-\ell})}
= \abs{1+e^{ik\tet_1}+z^\ell(1+z^{m-\ell})} \\
& = \abs{1-1+z^\ell(1+z^{m-\ell})}
= \abs{z^\ell(1+z^{m-\ell})} = \abs{1+z^{m-\ell}}
\end{align*}
which, by Lemma 6.1, is at most $\pi g/k$.

Now to prove $\mu(5)\le1+\pi/5$,
let $(0,a,b,c,d)\in\bN''_5$, and define
$\alf(z)=1+z^a+z^b+z^c+z^d$.
We must show that
$$
\abs{\alf(z)} \le 1 + \frac\pi5
$$
for some $z$ on the unit circle.
From the definition of $\bN''_5$, we have not only
$0<a<b<c<d$, but also $c\ge d-a$, which implies $c>d-b$.
We thus have the following inequalities, where both sides
are positive integers.
\begin{align*}
d & > b-a \\
d & > c-b \\
d & > c-a \\
c & > b-a \\
c & > d-b \\
c & \ge d-a
\end{align*}
We now define
\begin{align*}
g_1 & = \gcd(d,b-a) \qquad \mbox{so $g_1 \le b-a < d$} \\
g_2 & = \gcd(d,c-b) \qquad \mbox{so $g_2 \le c-b < d$} \\
g_3 & = \gcd(d,c-a) \qquad \mbox{so $g_3 \le c-a < d$} \\
g_4 & = \gcd(c,b-a) \qquad \mbox{so $g_4 \le b-a < c$} \\
g_5 & = \gcd(c,d-b) \qquad \mbox{so $g_5 \le d-b < c$} \\
g_6 & = \gcd(c,d-a) \qquad \mbox{so $g_6 \le d-a \le c$} \\
\end{align*}
Applying Lemma 6.2 to the case $(k,\ell,m)=(d,a,b)$,
we have
\begin{align*}
\abs{\alf(z)} & \le \abs{z^c} + \abs{1+z^d+z^a+z^b} \\
& = 1 + \abs{1+z^d+z^a+z^b} \le 1 + \frac{\pi g_1}{d} \qquad \mbox{for some $\abs{z}=1$}.
\end{align*}
Similarly, applying Lemma 6.2 to each of the cases
\begin{align*}
(k,\ell,m) & = (d,b,c) \\
(k,\ell,m) & = (d,a,c) \\
(k,\ell,m) & = (c,a,b) \\
(k,\ell,m) & = (c,b,d) \\
(k,\ell,m) & = (c,a,d)
\end{align*}
we can conclude that we have
\begin{align*}
\abs{\alf(z)} \le 1 + \frac{\pi g_2}{d} \qquad \mbox{for some $\abs{z}=1$} \\[1ex]
\abs{\alf(z)} \le 1 + \frac{\pi g_3}{d} \qquad \mbox{for some $\abs{z}=1$} \\[1ex]
\abs{\alf(z)} \le 1 + \frac{\pi g_4}{c} \qquad \mbox{for some $\abs{z}=1$} \\[1ex]
\abs{\alf(z)} \le 1 + \frac{\pi g_5}{c} \qquad \mbox{for some $\abs{z}=1$} \\[1ex]
\abs{\alf(z)} \le 1 + \frac{\pi g_6}{c} \qquad \mbox{for some $\abs{z}=1$}
\end{align*}
If any of the six numbers
\begin{equation} \label{sixthings}
\frac{g_1}{d}, \; \frac{g_2}{d}, \; \frac{g_3}{d}, \; \frac{g_4}{c}, \; \frac{g_5}{c}, \; \frac{g_6}{c}
\end{equation}
is less than or equal to $1/5$, we are done.
So for the remainder of the proof, we assume each of those six numbers
is greater than $1/5$.

Note that $d$ and $b-a$ are integer multiples of $g_1$, with $d > b-a$.
The condition $g_1/d > 1/5$ is equivalent to $d < 5g_1$, which implies
that $d \in \{2g_1,3g_1,4g_1\}$.
Then since $b-a$ is an integer multiple of $g_1$ smaller than $d$,
one of the following must be true:
\begin{align*}
\mbox{$d=2g_1$ and $b-a=g_1$} & \implies \frac12d = b-a \\[1ex]
\mbox{$d=3g_1$ and $b-a=g_1$} & \implies \frac13d = b-a \\[1ex]
\mbox{$d=3g_1$ and $b-a=2g_1$} & \implies \frac23d = b-a \\[1ex]
\mbox{$d=4g_1$ and $b-a=1g_1$} & \implies \frac14d = b-a \\[1ex]
\mbox{$d=4g_1$ and $b-a=3g_1$} & \implies \frac34d = b-a
\end{align*}
So we have $rd= b-a$, where
$$
r \in \Big\{ \frac14, \; \frac13, \; \frac12, \; \frac23, \; \frac34 \Big\}.
$$
We will now apply similar reasoning to the other five expressions in (\ref{sixthings}).
Note
\begin{gather*}
\mbox{$d$ and $c-b$ are integer multiples of $g_2$ with $d > c-b$, and $d<5g_2$} \\
\mbox{$d$ and $c-a$ are integer multiples of $g_3$ with $d > c-a$, and $d<5g_3$} \\
\mbox{$c$ and $b-a$ are integer multiples of $g_4$ with $c > b-a$, and $c<5g_4$} \\
\mbox{$c$ and $d-b$ are integer multiples of $g_5$ with $c > d-b$, and $c<5g_5$} \\
\mbox{$c$ and $d-a$ are integer multiples of $g_6$ with $c \ge d-a$, and $c<5g_6$}
\end{gather*}
Notice the slight difference in the last line. This means that our list of possibilities
for $c$ and $d-a$ will include $c=d-a=g_6$.

Using similar reasoning as before, we can conclude
\begin{gather*}
\mbox{$rd=b-a$, where $r \in \{\frac14, \; \frac13, \; \frac12, \; \frac23, \; \frac34\}$} \\
\mbox{$sd=c-b$, where $s \in \{\frac14, \; \frac13, \; \frac12, \; \frac23, \; \frac34\}$} \\
\mbox{$td=c-a$, where $t \in \{\frac14, \; \frac13, \; \frac12, \; \frac23, \; \frac34\}$} \\
\mbox{$uc=b-a$, where $u \in \{\frac14, \; \frac13, \; \frac12, \; \frac23, \; \frac34\}$} \\
\mbox{$vc=d-b$, where $v \in \{\frac14, \; \frac13, \; \frac12, \; \frac23, \; \frac34\}$} \\
\mbox{$wc=d-a$, where $w \in \{\frac14, \; \frac13, \; \frac12, \; \frac23, \; \frac34, \; 1\}$}
\end{gather*}
Thus there are finitely many possibilities for $r,s,t,u,v,w$.
Notice that this remains true if we change our goal from
showing $\mu(5)\le1+\pi/5$ to showing $\mu(5)\le1+\pi/m$.
In that case, the possible values for $r,s,t,u,v,w$ are
fractions between 0 and 1 whose denominators are less than $m$.

We claim that for each of the finitely many possible triples $r,s,u$,
there is at most one element of $\bN''_5$ satisfying the three conditions
\begin{align*}
rd & = b-a \\
sd & = c-b \\
uc & = b-a
\end{align*}
To verify this claim, we write these conditions as
\begin{align*}
1a - 1b + 0c + rd & = 0 \\
0a + 1b - 1c + sd & = 0 \\
1a - 1b + uc + 0d & = 0
\end{align*}
which we temporarily regard as a system of equations
in {\bf rational} unknowns $a,b,c,d$.
Writing as a matrix and row-reducing, we have
$$
\begin{bmatrix}
1 & -1 & 0 & r \\
0 & 1 & -1 & s \\
1 & -1 & u & 0
\end{bmatrix}
\rightarrow
\begin{bmatrix}
1 & -1 & 0 & r \\
0 & 1 & -1 & s \\
0 & 0 & u & -r
\end{bmatrix}
\rightarrow
\begin{bmatrix}
1 & -1 & 0 & r \\
0 & 1 & -1 & s \\
0 & 0 & 1 & -\frac{r}{u}
\end{bmatrix}
$$
Note that $u$ is always nonzero. Continuing, we have
$$
\begin{bmatrix}
1 & -1 & 0 & r \\
0 & 1 & -1 & s \\
0 & 0 & 1 & -\frac{r}{u}
\end{bmatrix}
\rightarrow
\begin{bmatrix}
1 & -1 & 0 & r \\
0 & 1 & 0 & s-\frac{r}{u} \\
0 & 0 & 1 & -\frac{r}{u}
\end{bmatrix}
\rightarrow
\begin{bmatrix}
1 & 0 & 0 & r+s-\frac{r}{u} \\
0 & 1 & 0 & s-\frac{r}{u} \\
0 & 0 & 1 & -\frac{r}{u}
\end{bmatrix}
$$
It follows that the system has an infinite one-parameter family
of {\bf rational} solutions, given by
\begin{align*}
a & = \Big(\frac{r}{u}-r-s\Big)q \\[1ex]
b & = \Big(\frac{r}{u}-s\Big)q \\[1ex]
c & = \Big(\frac{r}{u}\Big)q \\[1ex]
d & = q
\end{align*}
where $q$ can be any rational number.
In other words, the rational solutions $(a,b,c,d)$
are precisely the rational multiples of the 4-tuple
$$
\Big( \frac{r}{u}-r-s, \; \frac{r}{u}-s, \; \frac{r}{u}, \; 1 \Big)
$$
However, at most one rational multiple of a given rational 4-tuple
can consist of relatively prime nonnegative integers.
That is, there is at most one possible $(0,a,b,c,d)\in\bN''_5$
for each choice of $r,s,u$.

For the specific case of showing $\mu(5)\le1+\pi/5$, it is possible
to enumerate all the possibilities by hand. Recall that we have the
restrictions
\begin{gather*}
\mbox{$rd=b-a$, where $r \in \{\frac14, \; \frac13, \; \frac12, \; \frac23, \; \frac34\}$} \\
\mbox{$sd=c-b$, where $s \in \{\frac14, \; \frac13, \; \frac12, \; \frac23, \; \frac34\}$} \\
\mbox{$td=c-a$, where $t \in \{\frac14, \; \frac13, \; \frac12, \; \frac23, \; \frac34\}$} \\
\mbox{$uc=b-a$, where $u \in \{\frac14, \; \frac13, \; \frac12, \; \frac23, \; \frac34\}$} \\
\mbox{$vc=d-b$, where $v \in \{\frac14, \; \frac13, \; \frac12, \; \frac23, \; \frac34\}$} \\
\mbox{$wc=d-a$, where $w \in \{\frac14, \; \frac13, \; \frac12, \; \frac23, \; \frac34, \; 1\}$}
\end{gather*}
Notice that
$$
(r+s)d = rd + sd = (b-a) + (c-b) = c-a = td
$$
so $r+s=t$, and also notice that
$$
(u+v)c = uc + vc = (b-a) + (d-b) = d-a = wc
$$
so $u+v=w$.
The only permissible values of $r,s,t$ that satisfy this are
\begin{center}
\begin{tabular}{ccc}
$r$ & $s$ & $t$ \\ \hline
$1/4$ & $1/4$ & $1/2$ \\
$1/4$ & $1/2$ & $3/4$ \\
$1/3$ & $1/3$ & $2/3$ \\
$1/2$ & $1/4$ & $3/4$
\end{tabular}
\end{center}
and the only permissible values of $u,v,w$ that satisfy this are
\begin{center}
\begin{tabular}{ccc}
$u$ & $v$ & $w$ \\ \hline
$1/4$ & $1/4$ & $1/2$ \\
$1/4$ & $1/2$ & $3/4$ \\
$1/4$ & $3/4$ & $1$ \\
$1/3$ & $1/3$ & $2/3$ \\
$1/3$ & $2/3$ & $1$ \\
$1/2$ & $1/4$ & $3/4$ \\
$1/2$ & $1/2$ & $1$ \\
$2/3$ & $1/3$ & $1$ \\
$3/4$ & $1/4$ & $1$
\end{tabular}
\end{center}
Next, notice that we have $rd = b-a = uc$, which together with $c<d$
implies that $u>r$. Also notice that we have
\begin{gather*}
sd = c-b \implies b = c-sd \\
vc = d-b \implies b = d-vc \\
c-sd = d-vc \implies c+vc = d+sd \implies (1+v)c = (1+s)d
\end{gather*}
which together with $c<d$ implies $1+v > 1+s$, so $v>s$.
The only permissible values of $r,s,u,v$ satisfying both $u>r$ and $v>s$ are
\begin{center}
\begin{tabular}{cccc}
$r$ & $s$ & $u$ & $v$ \\ \hline
$1/4$ & $1/4$ & $1/3$ & $1/3$ \\
$1/4$ & $1/4$ & $1/3$ & $2/3$ \\
$1/4$ & $1/4$ & $1/2$ & $1/2$ \\
$1/4$ & $1/4$ & $2/3$ & $1/3$ \\
$1/4$ & $1/2$ & $1/3$ & $2/3$ \\
$1/3$ & $1/3$ & $1/2$ & $1/2$ \\
$1/2$ & $1/4$ & $2/3$ & $1/3$ \\
\end{tabular}
\end{center}
We now use the conditions $rd=uc$ and $(1+v)c = (1+s)d$
to conclude
$$
\frac{r}{u} = \frac{c}{d} = \frac{1+s}{1+v}
\implies r(1+v) = u(1+s)
$$
which only some of our possible values of $r,s,u,v$
will satisfy.
\begin{center}
\begin{tabular}{cccccccc}
$r$ & $s$ & $u$ & $v$ & $1+s$ & $1+v$ & $r(1+v)$ & $u(1+s)$ \\ \hline
$1/4$ & $1/4$ & $1/3$ & $1/3$ & $5/4$ & $4/3$ & $1/3$ & $5/12$ \\
$1/4$ & $1/4$ & $1/3$ & $2/3$ & $5/4$ & $5/3$ & $5/12$ & $5/12$ \\
$1/4$ & $1/4$ & $1/2$ & $1/2$ & $5/4$ & $3/2$ & $3/8$ & $5/8$ \\
$1/4$ & $1/4$ & $2/3$ & $1/3$ & $5/4$ & $4/3$ & $1/3$ & $5/6$ \\
$1/4$ & $1/2$ & $1/3$ & $2/3$ & $3/2$ & $5/3$ & $5/12$ & $1/2$ \\
$1/3$ & $1/3$ & $1/2$ & $1/2$ & $4/3$ & $3/2$ & $1/2$ & $2/3$ \\
$1/2$ & $1/4$ & $2/3$ & $1/3$ & $5/4$ & $4/3$ & $2/3$ & $5/6$
\end{tabular}
\end{center}
We see that only one of our possibilities satisfies $r(1+v)=u(1+s)$.
We conclude that $(r,s,u,v)=(\frac14,\frac14,\frac13,\frac23)$.
We then have
$$
\Big(\frac{r}{u}-r-s,\;\frac{r}{u}-s,\;\frac{r}{u},\;1\Big)
=\Big(\frac14,\;\frac12,\;\frac34,\;1\Big)
$$
implying $(a,b,c,d)=(1,2,3,4)$. But the polynomial
$f(z)=1+z+z^2+z^3+z^4$ certainly satisfies $\abs{f(z)}\le1+\pi/5$
for some $\abs{z}=1$, because it has zeros at the nontrivial
fifth roots of unity.
This completes the proof that $\mu(5)\le1+\pi/5$.

The above argument can be modified to show that $\mu(5)\le1+\pi/m$
for some other positive integers $m$, but doing so by hand is cumbersome
and computer assistance is helpful. We will give a rough outline
of an argument that $\mu(5)\le1+\pi/6$.

This time, the argument involves finding $r,s,t,u,v,w$ that satisfy
\begin{align*}
rd & = b-a \\
sd & = c-b \\
td & = c-a \\
uc & = b-a \\
vc & = d-b \\
wc & = d-a
\end{align*}
where $0<r,s,t,u,v<1$, $0<w\le1$, and $r,s,t,u,v,w$ are fractions
with denominators strictly less than 6.
They must further satisfy
$$
r+s=t, \qquad u+v=w, \qquad r<u, \qquad s<v, \qquad r(1+v)=u(1+s).
$$
A finite search (aided by computer) reveals that the only eligible
values of $r,s,u,v$ are those shown in this table.
(The values of $1+s$ and $1+v$ are included for convenience.)
\begin{center}
\begin{tabular}{cccccc}
$r$ & $s$ & $u$ & $v$ & $1+s$ & $1+v$ \\ \hline
$1/3$ & $1/3$ & $2/5$ & $3/5$ & $4/3$ & $8/5$ \\
$1/4$ & $1/4$ & $1/3$ & $2/3$ & $5/4$ & $5/3$ \\
$1/5$ & $1/5$ & $1/4$ & $1/2$ & $6/5$ & $3/2$ \\
$1/5$ & $2/5$ & $1/4$ & $3/4$ & $7/5$ & $7/4$ \\
$2/5$ & $1/5$ & $1/2$ & $1/2$ & $6/5$ & $3/2$ \\
$3/5$ & $1/5$ & $2/3$ & $1/3$ & $6/5$ & $4/3$
\end{tabular}
\end{center}
As before, each eligible $(r,s,u,v)$ gives us a 4-tuple
$$
\Big(\frac{r}{u}-r-s,\;\frac{r}{u}-s,\;\frac{r}{u},\;1\Big)
$$
of which precisely one integer multiple is an eligible $(a,b,c,d)$.
This gives us six possibilities.
\begin{center}
\begin{tabular}{cccc}
$\frac{r}{u}-r-s$ & $\frac{r}{u}-s$ & $\frac{r}{u}$ & $(a,b,c,d)$ \\ \hline
$1/6$ & $3/6$ & $5/6$ & $(1,3,5,6)$ \\
$1/4$ & $2/4$ & $3/4$ & $(1,2,3,4)$ \\
$2/5$ & $3/5$ & $4/5$ & $(2,3,4,5)$ \\
$1/5$ & $2/5$ & $4/5$ & $(1,2,4,5)$ \\
$1/5$ & $3/5$ & $4/5$ & $(1,3,4,5)$ \\
$1/10$ & $7/10$ & $9/10$ & $(1,7,9,10)$ \\
\end{tabular}
\end{center}
One can verify that for each of these six possibilities for $(a,b,c,d)$,
we have $\abs{1+z^a+z^b+z^c+z^d}\le1+\pi/6$ for some $\abs{z}=1$.
This completes the sketch of the proof that $\mu(5)\le1+\pi/6$.

In conclusion, note that the contributions of Campbell et al.
and Goddard appeared gradually. The results in this paper can be
regarded as extensions of that work. It appears that evaluating
$\lam(n)$ or $\mu(n)$ in finitely many steps is a genuinely subtle problem.

\bibliographystyle{amsplain}

\begin{thebibliography}{10}

\bibitem{Boy}
D.W. Boyd, {\it Large Newman polynomials}, in {\it Diophantine analysis (Kensington, 1985)},
159--170, London Math.\ Soc.\ Lecture Note Ser., 109,
Cambridge Univ.\ Press, Cambridge, 1986.

\bibitem{CFF}
D.M. Campbell, H.R.P. Ferguson \& R.W. Forcade,
{\it Newman polynomials on $\abs{z}=1$},
Indiana Univ.\ Math.\ J. {\bf 32} (1983), 517--525.

\bibitem{Cho52}
S. Chowla,
{\it The Riemann zeta and allied functions},
Bull.\ Amer.\ Math.\ Soc.\ {\bf 58} (1952), 287--305.

\bibitem{Cho65}
S. Chowla,
{\it Some applications of a method of A. Selberg},
J. Reine Angew.\ Math.\ {\bf 217} (1965), 128--132.

\bibitem{Coh}
P.J. Cohen,
{\it On a conjecture of Littlewood and idempotent measures},
Amer.\ J. Math.\ {\bf 82} (1960), 191--212.

\bibitem{God}
B. Goddard,
{\it Finite exponential series and Newman polynomials},
Proc.\ Amer.\ Math.\ Soc.\ {\bf 116} (1992), 313--320.

\bibitem{Mer}
I. Mercer,
{\it Newman polynomials not vanishing on the unit circle},
Integers {\bf 12} (2012), A67, 7 pp.

\bibitem{Ruz}
I.Z. Ruzsa,
{\it Negative values of cosine sums},
Acta Arith. {\bf 111} (2004), 179--186.

\bibitem{UU}
M. Uchiyama (n\'{e}e Katayama) \& S. Uchiyama,
{\it On the cosine problem},
Proc.\ Japan Acad.\ {\bf 36} (1960), 475--479.

\end{thebibliography}

\end{document}